\documentclass[12pt]{article}
\usepackage{amssymb}

\newcommand{\be}{\begin{equation}}
\newcommand{\ee}{\end{equation}}
\newcommand{\bea}{\begin{eqnarray}}
\newcommand{\eea}{\end{eqnarray}}
\newcommand{\barray}{\begin{array}}
\newcommand{\earray}{\end{array}}

\newcommand{\nn}{\nonumber}
\newcommand{\bitem}{\begin{itemize}}
\newcommand{\eitem}{\end{itemize}}
\newtheorem{teo}{Theorem}[section]
\newcommand{\bt}{\begin{teo}}
\newcommand{\et}{\end{teo}}
\newtheorem{Def}{Definition}[section]
\newcommand{\bd}{\begin{Def}}
\newcommand{\ed}{\end{Def}}
\newtheorem{lem}{Lemma}[section]
\newcommand{\bl}{\begin{lem}}
\newcommand{\el}{\end{lem}}
\newtheorem{prop}{Proposition}[section]
\newcommand{\bp}{\begin{prop}}
\newcommand{\ep}{\end{prop}}
\newtheorem{cor}{Corollary}[section]
\newcommand{\bc}{\begin{cor}}
\newcommand{\ec}{\end{cor}}
\newtheorem{ex}{Example}[section]
\newcommand{\bex}{\begin{ex}}
\newcommand{\eex}{\end{ex}}
\newtheorem{rem}{Remark}[section]
\newcommand{\br}{\begin{rem}}
\newcommand{\er}{\end{rem}}

\begin{document}

\begin{center}
{\Large \textbf{On commutative subalgebras of the Weyl algebra\\
that are related to commuting operators \\of arbitrary rank 
and genus\footnotetext[1]{This research was supported by 
the Russian Foundation for Basic Research (grant nos. 11-01-00197
and 11-01-12067-ofi-m-2011) and by the ``Leading Scientific Schools'' programme (grant
no. NSh-4995-2012.1).}}}
\end{center}

\medskip
\smallskip

\begin{center}
{\large \bf {O. I. Mokhov}}
\end{center}

\bigskip
\bigskip
\medskip

In this note we construct examples of commuting ordinary scalar differential
operators with polynomial coefficients that are related to a spectral curve of an arbitrary genus
$g$ and to an arbitrary even rank $r = 2k$, and also to an arbitrary rank of the form $r = 3k$, of
the vector bundle of common eigenfunctions of the commuting operators over the spectral curve.
At present we know no explicit examples of commuting ordinary scalar differential
operators that are related to a spectral curve of genus
$g > 1$ for rank of the form $r = 6s \pm 1$, $s \geq 1$. For all other values of
genus $g$ and rank $r$, explicit examples of commuting
operators even with polynomial coefficients are constructed.
We conjecture that there exist commuting ordinary scalar differential
operators with polynomial coefficients that are related to a spectral curve of an arbitrary genus
$g > 1$ also for an arbitrary rank of the form $r = 6s \pm 1$, $s \geq 1$, but such examples are not known for now, -- this is a very interesting problem.
It is well known that commuting ordinary scalar differential
operators with polynomial coefficients give commutative subalgebras of the Weyl algebra
$W_1$, i.e., the algebra with two generators $p$ and $q$ and the relation
$[p, q] = 1$.
By the Burchnall--Chaundy lemma [1] any pair of commuting ordinary scalar differential
operators $L$ and $M$ is connected by a certain polynomial relation
$Q (L, M) = 0$ given by the spectral curve $Q(\lambda, \mu) = 0$ of the pair of commuting
operators:
$L \psi = \lambda \psi$, $M \psi = \mu \psi$, and common eigenfunctions of the commuting operators
define an $r$-dimensional vector bundle over the spectral curve
(the dimension $r$ of
the vector bundle of common eigenfunctions of a pair of commuting operators over the spectral curve
is called {\it the rank of the
commuting pair of operators};
the rank of any pair of commuting operators is a common divisor of the orders of these
commuting operators).
The first examples of commuting ordinary scalar differential
operators of the nontrivial ranks 2 and 3 and the nontrivial genus $g = 1$
were constructed by Dixmier [2] for the nonsingular elliptic spectral curve
$\mu^2 = \lambda^3 - \alpha$,
where $\alpha$ is an arbitrary nonzero constant:
\be
L = \left (  \left ({d \over dx}\right )^2 + x^3 + \alpha \right )^2 + 2x,
\ee
\be
M = \left (  \left ({d \over dx}\right )^2 + x^3 + \alpha \right)^3 + 3x \left ( {d \over dx} \right)^2 + 3 {d \over dx}
+ 3x (x^3 + \alpha),
\ee
where $L$ and $M$ is the commuting pair of the Dixmier operators of rank 2, genus 1, $M^2 = L^3 - \alpha$, $[L, M] = 0$;
\be
L = \left (  \left ({d \over dx}\right )^3 + x^2 + \alpha \right )^2 + 2 {d \over dx},
\ee
\be
M = \left (  \left ({d \over dx}\right )^3 + x^2 + \alpha \right)^3 + 3 \left ( {d \over dx} \right)^4 + 3 (x^2 + \alpha) {d \over dx} + 3x,
\ee
where $L$ and $M$ is the commuting pair of the Dixmier operators of rank 3, genus 1, $M^2 = L^3 - \alpha$, $[L, M] = 0$.
These remarkable examples were found by Dixmier as commutative subalgebras of the Weyl algebra $W_1$ [2].
The general classification of commuting ordinary scalar differential
operators of nontrivial ranks $r > 1$ was obtained by Krichever [3]. The general form of commuting
ordinary scalar differential operators of rank 2 for an arbitrary elliptic spectral curve was found
by Krichever and Novikov [4]. The description of commuting ordinary scalar differential
operators with polynomial coefficients (commutative subalgebras of the Weyl algebra $W_1$) is a separate nontrivial
problem, and this problem is not solved completely even for the Krichever--Novikov commuting operators of rank 2, genus 1, which are rationally parametrized by one arbitrary function (this problem was considered and studied in [5], [6]). We will also devote a separate paper to this very interesting question. The functional parameter corresponding to the Dixmier example of rank 2, genus 1 among all
the Krichever--Novikov commuting operators of rank 2, genus 1 was found by Grinevich [7]. The general form of commuting
ordinary scalar differential operators of rank 3 for an arbitrary elliptic spectral curve (the general commuting operators of rank 3, genus 1 are parametrized
by two arbitrary functions) was found by Mokhov [8], [9], where were also found the functional parameters corresponding to the Dixmier example of rank 3, genus 1 among all the commuting operators of rank 3, genus 1. Moreover, examples of commuting ordinary scalar differential
operators of genus 1 with polynomial coefficients were constructed for any rank $r$ (some commutative subalgebras of the Weyl algebra $W_1$) (see also [10]). We note that the coefficients of one from a pair of commuting operators are always expressed polynomially in terms of the coefficients of the second operator of the commuting pair and their derivatives.
Recently Mironov [11] (see also earlier papers [12]--[14]) constructed for any genus $g > 1$ remarkable
examples of commuting ordinary scalar differential
operators of ranks 2 and 3 with polynomial coefficients that generalize naturally the Dixmier examples of ranks 2 and 3, genus 1:
\be
L = \left (  \left ({d \over dx}\right )^2 + x^3 + \alpha \right )^2 + g (g + 1) x,
\ee
\be
M^2 = L^{2g + 1} + a_{2g} L^{2g} + \ldots + a_1 L + a_0,
\ee
where $a_i$ are some constants, $\alpha$ is an arbitrary nonzero constant, $L$ and $M$ are the Mironov commuting operators of rank 2, genus $g$ (the orders of the operators $L$ and $M$ are 4 and $4g + 2$, respectively), the coefficients of the operator $M$ are expressed polynomially in terms of the coefficients of the operator $L$ and their derivatives, $[L, M] = 0$;
\be
L = \left (  \left ({d \over dx}\right )^3 + x^2 + \alpha \right )^2 + g (g + 1) {d \over dx},
\ee
\be
M^2 = L^{2g + 1} + a_{2g} L^{2g} + \ldots + a_1 L + a_0,
\ee
where $a_i$ are some constants, $\alpha$ is an arbitrary nonzero constant, $L$ and $M$ are the Mironov commuting operators of rank 3, genus $g$ (the orders of the operators $L$ and $M$ are 6 and $6g + 3$, respectively), the coefficients of the operator $M$ are expressed polynomially in terms of the coefficients of the operator $L$ and their derivatives, $[L, M] = 0$.

Let us construct examples of commuting ordinary scalar differential
operators with polynomial coefficients that are related to a spectral curve of an arbitrary genus
$g$ for an arbitrary even rank $r = 2k$, $k > 1$.

{\bf Theorem 1.} {\it
The operators $L$ and $M$ of orders $4k$ and $4kg + 2k$, respectively,
\be
L = \left (  \left ({d \over dx}\right )^{2k} - 2x \left ({d \over dx}\right )^{k} - k \left ({d \over dx}\right )^{k - 1} + \left ({d \over dx}\right )^3 + x^2 + \alpha \right )^2 + g (g + 1) {d \over dx},
\ee
\be
M^2 = L^{2g + 1} + a_{2g} L^{2g} + \ldots + a_1 L + a_0,
\ee
where $a_i$ are some constants, $\alpha$ is an arbitrary nonzero constant, are commuting operators of rank $r = 2k$, $k > 1$, genus $g$, $[L, M] = 0$, the coefficients of the operator $M$ are expressed polynomially in terms of the coefficients of the operator $L$ and their derivatives. For $k >2$ the commuting operators $L$ and $M$ have the standard {\rm(}canonical{\rm)} form.}

Let us also construct examples of commuting ordinary scalar differential
operators with polynomial coefficients that are related to a spectral curve of an arbitrary genus
$g$ for an arbitrary rank of the form $r = 3k$, $k \geq 1$.

{\bf Theorem 2.} {\it
The operators $L$ and $M$ of orders $6k$ and $6kg + 3k$, respectively,
\bea
&& L = \left (  \left ({d \over dx}\right )^{3k} - 3x \left ({d \over dx}\right )^{2k} - 3k \left ({d \over dx}\right )^{2k - 1} + 3 x^2 \left ({d \over dx}\right )^k + 3 k x \left ({d \over dx}\right )^{k - 1} + \right. \nn\\ && \left. + k (k - 1) \left ({d \over dx}\right )^{k - 2} + \left ({d \over dx}\right )^2 - x^3 + \alpha \right )^2 - g (g + 1) x,
\eea
\be
M^2 = L^{2g + 1} + a_{2g} L^{2g} + \ldots + a_1 L + a_0,
\ee
where $a_i$ are some constants, $\alpha$ is an arbitrary nonzero constant, are commuting operators of rank $3k$, $k \geq 1$, genus $g$, $[L, M] = 0$, the coefficients of the operator $M$ are expressed polynomially in terms of the coefficients of the operator $L$ and their derivatives. For $k >1$ the commuting operators $L$ and $M$ have the standard {\rm(}canonical{\rm)} form.}

\smallskip

\begin{center}
\bf {References}
\end{center}

\smallskip

[1] J.L. Burchnall, I.W. Chaundy. Commutative ordinary differential operators. {\it Proc. London Math. Society}, ser. 2, {\bf 21} (1923), 420--440.

%[2] J.L. Burchnall, I.W. Chaundy, {\it Royal Soc. London}, ser. A, {\bf 118} (1928), 557--583.
%[3] J.L. Burchnall, I.W. Chaundy, {\it Royal Soc. London}, ser. A, {\bf 134} (1931), 471--485.

[2] J. Dixmier. Sur les alg\`ebres de Weyl. {\it Bull. Soc. Math. France}, {\bf 96} (1968), 209--242.

[3] I.M. Krichever. Commutative rings of ordinary linear differential operators. {\it Funktsional. Anal. i Prilozhen.}, {\bf 12}:3 (1978), 20--31 (in Russian); English translation in {\it Functional Anal. Appl.}, {\bf 12}:3 (1978), 175--185.

[4] I.M. Krichever, S.P. Novikov. Holomorphic bundles and nonlinear equations. Finite-zone solutions of rank 2. {\it Dokl. Akad. Nauk SSSR}, {\bf 247}:1 (1979), 33--36 (in Russian); English translation in {\it Soviet Math. Dokl.}.

[5] O.I.Mokhov. On commutative subalgebras of Weyl algebra, which are associated with an elliptic curve.
International Conference on Algebra in Memory of A.I.Shirshov (1921--1981). Barnaul, USSR, 20--25 August 1991.
Reports on theory of rings, algebras and modules. 1991. P. 85.

[6] O.I.Mokhov. On the commutative subalgebras of Weyl algebra, which are generated by the Chebyshev polynomials.
Third International Conference on Algebra in Memory of M.I.Kargapolov (1928--1976). Krasnoyarsk, Russia, 23--28 August 1993. Krasnoyarsk: Inoprof, 1993. P. 421.

[7] P.G. Grinevich. Rational solutions for the equation of commutation of differential operators. {\it Funktsional. Anal. i Prilozhen.}, {\bf 16}:1 (1982), 19--24 (in Russian); English translation in {\it Functional Anal. Appl.}, {\bf 16}:1 (1982), 15--19.

[8] O.I. Mokhov. Commuting ordinary differential operators of rank 3 corresponding to an elliptic curve. {\it Uspekhi Matem. Nauk}, {\bf 37}:4 (1982), 169--170 (in Russian); English translation in {\it Russian Math. Surveys}, {\bf 37}:4 (1982), 129--130.

[9] O.I. Mokhov. Commuting differential operators of rank 3, and nonlinear differential equations. {\it Izvestiya AN SSSR, ser. matem.}, {\bf 53}:6 (1989), 1291--1315 (in Russian); English translation in {\it Math. USSR, Izvestiya}, {\bf 35}:3 (1990), 629--655.

[10] P. Dehornoy. Op\'erateurs diff\'erentiels et courbes elliptiques. {\it Compositio Math.}, {\bf 43}:1 (1981), 71--99.

[11] A.E. Mironov. Self-adjoint commuting differential operators and commutative subalgebras of the Weyl algebra.
arXiv: 1107.3356.

[12]  A.E. Mironov. On commuting differential operators of rank 2. {\it Sib. Elektron. Matem. Izvestiya [Siberian Electronic Mathematical Reports]}, {\bf 6} (2009), 533--536 (in Russian).

[13] A.E. Mironov. Commuting rank 2 differential operators corresponding to a curve of genus 2. {\it Funktsional. Anal. i Prilozhen.}, {\bf 39}:3 (2005), 91--94 (in Russian); English translation in
{\it Functional Anal. Appl.}, {\bf 39}:3 (2005), 240--243.

[14] A.E.Mironov. A ring of commuting differential operators of rank 2 corresponding to a curve of genus 2.
{\it Matem. Sbornik}, {\bf 195}:5 (2004), 103--114 (in Russian); English translation in
{\it Sbornik: Mathematics}, {\bf 195}:5 (2004), 711--722.

\begin{flushleft}
%{\bf O.I. Mokhov}\\
Department of Geometry and Topology,\\
Faculty of Mechanics and Mathematics,\\
Lomonosov Moscow State University\\
Moscow, 119991 Russia\\
{\it E-mail\,}: mokhov@mi.ras.ru; mokhov@landau.ac.ru; mokhov@bk.ru\\
\end{flushleft}
\end{document}